\documentclass[A4paper, 12pt, french]{article}
\usepackage[utf8]{inputenc,vietnam}
\usepackage{cite}

\usepackage[portrait, top=2cm, bottom=2cm, left=2cm, right=2cm] {geometry}
\usepackage[french]{babel}
\usepackage{url}
\usepackage{array,epsfig}
\usepackage{amsmath}
\usepackage{amsfonts}
\usepackage{amssymb}
\usepackage{amsxtra}
\usepackage[thmmarks, amsmath, thref, amsthm]{ntheorem}
\usepackage{latexsym}
\usepackage{dsfont}
\usepackage{mathrsfs}
\usepackage{color}
\usepackage[all]{xy}
\usepackage{tensor}
\usepackage[pagebackref=true]{hyperref}

\usepackage{graphicx}
\usepackage{mathtools}
\usepackage{caption}
\usepackage[OT2, T1]{fontenc}
\usepackage{accents}
\usepackage[nottoc]{tocbibind}
\usepackage[page,toc,titletoc,title]{appendix}
\usepackage{sectsty}

\theoremstyle{plain}
\newtheorem{thm}{Th\'eor\`eme}
\newtheorem{lemm}[thm]{Lemme}
\newtheorem{exam}[thm]{Exemple}
\newtheorem{prop}[thm]{Proposition}
\newtheorem{coro}[thm]{Corollaire}

\newcommand{\Br}{\operatorname{Br}}
\newcommand{\Pic}{\operatorname{Pic}}
\newcommand{\Gal}{\operatorname{Gal}}
\newcommand{\SL}{\operatorname{SL}}
\newcommand{\Hom}{\operatorname{Hom}}
\newcommand{\Img}{\operatorname{Im}}
\newcommand{\Ker}{\operatorname{Ker}}
\newcommand{\res}{\operatorname{res}}
\newcommand{\Spec}{\operatorname{Spec}}
\newcommand{\Ext}{\operatorname{Ext}}
\newcommand{\Aut}{\operatorname{Aut}}

\newcommand{\PGCD}{\operatorname{PGCD}}
\newcommand{\ord}{\operatorname{ord}}
\newcommand{\et}{\textnormal{\'et}}
\newcommand{\ab}{\operatorname{ab}}
\newcommand{\grp}{\operatorname{grp}}
\newcommand{\type}{\texttt{type}}
\newcommand{\lien}{\texttt{lien}}

\renewcommand{\H}{\mathrm{H}}
\newcommand{\R}{\mathrm{R}}
\newcommand{\C}{\mathrm{C}}

\newcommand{\Cbb}{\mathbb{C}}
\newcommand{\Gbb}{\mathbb{G}}
\newcommand{\Pbb}{\mathbb{P}}
\newcommand{\Qbb}{\mathbb{Q}}
\newcommand{\Zbb}{\mathbb{Z}}

\newcommand{\Acal}{\mathcal{A}}
\newcommand{\Gcal}{\mathcal{G}}

\newcommand{\gfrak}{\mathfrak{g}}

\newcommand{\del}{\partial}

\newcommand{\ol}[1]{\overline{#1}}
\newcommand{\tuple}[1]{\left(#1\right)}
\newcommand{\pair}[1]{\left<#1\right>}
\newcommand{\ps}[1]{(\!(#1)\!)}
\newcommand{\pssq}[1]{[\![#1]\!]}

\DeclareSymbolFont{cyrletters}{OT2}{wncyr}{m}{n}
\DeclareMathSymbol{\Sha}{\mathalpha}{cyrletters}{"58}

\title{Groupes de Brauer alg\'ebriques modulo les constants d'espaces homog\`enes et leurs compactifications}
\author{Nguy$\tilde{\hat{\text{e}}}$n M$\d{\text{a}}$nh Linh}
\date{\today}

\begin{document}
	\maketitle
	\begin{abstract} 
  Soit $X$ une vari\'et\'e lisse, g\'eom\'etriquement int\`egre, sans fonctions inversibles non constantes sur un corps $K$. Alors le quotient du groupe Brauer \og alg\'ebrique \fg{} de $X$ par $\Br K$ s'injecte dans $\H^1(K,\Pic{\ol{X}})$. Nous montrons que cette inclusion n'est pas toujours un isomorphisme m\^eme dans le cas o\`u $X$ est un espace homog\`ene d'un groupe alg\'ebrique lin\'eaire connexe sur $K$. Un r\'esultat similaire pour les compactifications lisses de $X$ est aussi donn\'e.
\end{abstract}
    \bigskip

\section{Contexte du probl\`eme} 

Fixons quelques notations. Si $K$ est un corps, $\ol{K}$ d\'esigne une cl\^oture {\it s\'eparable} fix\'ee de $K$ et $\Gamma_K:=\Gal(\ol{K}/K)$. Pour $r \ge 0$, on note $\Qbb/\Zbb(r):=\varinjlim_n \mu_n^{\otimes r}$, et $\mu:=\Qbb/\Zbb(1)$ le sous-module de torsion de $\ol{K}^\times$. Si $A$ est un $\Gamma_K$-module, on note $A' = \Hom(A,\ol{K}^\times)$ son dual de Cartier. On d\'efinit
	$$\Sha^1_\omega(K,A):=\Ker\tuple{\H^1(K,A) \to \prod_{\sigma \in \Gamma_K} \H^1(\ol{\pair{\sigma}},A)}.$$
Si $L$ est une extension finie galoisienne de $K$ d\'eployant $A$, alors $\Sha^1_\omega(L,A) = 0$, et donc la suite exacte d'inflation-restriction donne
	\begin{equation} \label{eqSha1Omega}
		\Sha^1_\omega(K,A) = \Sha^1_\omega(\Gal(L/K),A) := \Ker\tuple{\H^1(\Gal(L/K),A) \to \prod_{g \in \Gal(L/K)} \H^1(\pair{g},A)}.
	\end{equation}

Soit $X$ une vari\'et\'e lisse et g\'eom\'etriquement int\`egre sur $K$. On note $\ol{X}:=X \times_K \ol{K}$. Par convention, le groupe de Brauer de $X$ est toujours le groupe de Brauer-Grothendieck $\Br X:=\H^2_{\et}(X,\Gbb_m)$. Le groupe de Brauer {\it alg\'ebrique} de $X$ est $\Br_1 X := \Ker(\Br X \to \Br \ol{X})$. La fl\`eche naturelle $\Br K \to \Br X$ se factorise par $\Br_1 X$ (puisque $\Br \ol{K} = 0$), et son image est le sous-groupe des \'el\'ements {\it constants} de $\Br_1 X$. On dispose d'une suite exacte
	\begin{equation} \label{eqExactSequenceX}
		0 \to \Br_1 X / \Br K \to \H^1(K,\Pic\ol{X}) \to \Ker(\H^3(K,\Gbb_m) \to \H^3_\et(X,\Gbb_m))
	\end{equation}
tir\'ee de la suite spectrale de Hochschild-Serre $\H^p(K,\H^q_\et(\ol{X},\Gbb_m)) \Rightarrow \H^{p+q}_\et(X,\Gbb_m)$. \`A la suite du travail de Borovoi, Demarche et Harari \cite[Th\'eor\`eme 8.1]{BDH}, on sait que dans le cas o\`u $K$ est de caract\'eristique nulle et o\`u $X$ est un espace homog\`ene d'un $K$-groupe lin\'eaire semi-simple simplement connexe, \`a stabilisateur g\'eom\'etrique r\'eductif, \eqref{eqExactSequenceX} devient
\begin{equation} \label{eqExactSequenceXc}
	0 \to \Br_1 X^c / \Br K \to \Sha^1_\omega(K,\Pic\ol{X}) \to \Ker(\H^3(K,\Gbb_m) \to \H^3_\et(X,\Gbb_m))
\end{equation}
o\`u $X^c$ d\'esigne une compatification lisse de $X$.

Lorsque $\H^3(K,\Gbb_m) = 0$ (par exemple si $K$ est un corps local, un corps global, ou encore le corps de fonction d'une courbe sur l'un de ces corps) ou si $X(K) \neq \varnothing$, la suite exacte \eqref{eqExactSequenceX} (resp. \eqref{eqExactSequenceXc}) donne un isomorphisme $\Br_1 X / \Br K \simeq \H^1(K,\Pic\ol{X})$ (resp. $\Br_1 X^c / \Br K \simeq \Sha^1_\omega(K,\Pic\ol{X})$). On est int\'eress\'e par la question suivante, soulev\'ee par Harari: Ces isomorphismes valent-t-ils sur un corps quelconque ? Le contexte de cette question est comme suit. Soit $k$ un corps de nombres et soit $f: X \to B$ un morphisme projectif, dominant de $k$-vari\'et\'es lisses, on veut \'etudier l'obstruction de Brauer-Manin pour l'espace total $X$, sachant celle pour la base $B$ et pour les fibres au-dessus d'\og assez \fg{} $k$-points de $B$. Il convient de se d\'eplacer entre les groupes de Brauer de ces fibres. Par exemple, cela est possible si des \og fl\`eches de sp\'ecialisations \fg{} $\Br_1 X_\eta / \Br k(\eta) \to \Br_1 X_b / \Br k$ sont des isomorphismes pour beaucoup de points $b \in B(k)$ (o\`u $\eta$ d\'esigne le point g\'en\'erique de $B$) - cette m\'ethode est le \og lemme formel \fg{} \cite[Corollaire 2.6.1]{harari1994fibration}. Pour l'argument de ce type, il est crucial que les \'elements de $\H^1(k(\eta),\Pic\ol{X_\eta})$ (resp. de $\Sha^1_\omega(k(\eta),\Pic\ol{X_\eta})$) proviennent de $\Br_1 X_\eta / \Br k(\eta)$ (resp. de $\Br_1 X_\eta^c / \Br k(\eta)$). Bien entendu, dans le cas o\`u $B$ est une courbe (sur le corps de nombres $k$) ou $f$ admet une section, ce probl\`eme dispara\^it. Harari ont trait\'e le cas o\`u $B = \Pbb^n_k$ - mais en utilisant un argument de r\'ecurrence qui s'appuie fortement sur la nature des espaces projectifs (qui ne marche pas m\^eme dans le cas o\`u $B$ est $k$-rationnelle).

Nous allons donner dans ce texte une r\'eponse n\'egative aux deux parties de la question mention\'ee ci-dessus (le corollaire \ref{corDifferential} et le th\'eor\`eme \ref{thmMain}). Pour le groupe $\Br_1 X / \Br K$, c'est une cons\'equence simple du th\'eor\`eme de Rost-Voevodsky (anciennement connu sous le nom de la conjecture de Bloch-Kato) d\`es que la fl\`eche $\H^1(K,\Pic\ol{X}) \to \H^3(K,\Gbb_m)$ de \eqref{eqExactSequenceX} est explicit\'ee dans le \S \ref{section1}. Pour le contre-exemple concernant le groupe $\Br_1 X^c / \Br K$, il s'agit de la dualit\'e sur des corps 2-locaux.

\section{Une diff\'erentielle de la suite spectrale de Hochschild-Serre} \label{section1}

Soient $K$ un corps, $G$ un $K$-groupe lin\'eaire semi-simple simplement connexe, et $X$ un espace homog\`ene de $G$. On ne suppose pas que $X$ poss\`ede un $K$-point. Notons $\ol{H}$ le stabilisateur d'un $\ol{K}$-point de $X$, qu'on suppose r\'eductif. Si $\ol{H}$ n'est pas ab\'elien, il n'est pas n\'ecessairement d\'efini sur $K$. Cependant, on peut toujours associe \`a $X$ le $K$-{\it lien de Springer} $L_X$ (dont le $\ol{K}$-groupe sous-jacent est $\ol{H}$), l'ensemble de 2-cohomologie non ab\'elien $\H^2(K,L_X)$, ainsi que la {\it classe de Springer} $\eta_X \in \H^2(K,L_X)$ (voir \cite[\S 1]{FSS} pour leur d\'efinition). La neutralit\'e de $\eta_X$ est une obstruction \`a l'existence d'espaces homog\`enes principaux de $G$ dominant $X$.

Il y a une action naturelle de $\Gamma_K$ sur l'ab\'elianis\'e $\ol{H}^{\ab} = \ol{H}/[\ol{H},\ol{H}]$. On note $H^{\ab}$ la $K$-forme correspondante, qui est un $K$-groupe de type multiplicatif ($\ol{H}$ \'etant r\'eductif). En particulier, $\H^0_\et(\ol{X},\ol{H}^{\ab}) = \ol{H}^{\ab}$. En effet, $\ol{H}^{\ab} = F \times \Gbb_m^r$, o\`u $F$ est certain $\ol{K}$-groupe \'etale et $r \ge 0$ est un entier. D'une part, $\H^0_\et(\ol{X},F) = F$ comme $\ol{X}$ est connexe. D'autre part $\ol{K}[G]^\times = \ol{K}^\times$ par le lemme de Rosenlicht \cite[Proposition 3]{Rosenlicht} ($G$ \'etant semisimple), {\it a fortiori} $\ol{K}[X]^\times = \ol{K}^\times$, d'o\`u $\H^0_\et(\ol{X},\Gbb_m^r) = (\ol{K}^\times)^r$.

Notons de plus que le dual de Cartier $H':=(H^{\ab})'$ est un $\Gamma_K$-module discret qui est de type fini en tant que groupe ab\'elien. La projection $\ol{H} \to \ol{H}^{\ab}$ induit un morphisme de $L_X \to \lien(H^{\ab})$ de $K$-liens alg\'ebriques. Par surjectivit\'e, ce morphisme induit une {\it application} $\H^2(K,L_X) \to \H^2(K,H^{\ab})$ ({\it a priori}, un morphisme de $K$-liens induit seulement une {\it relation} entre les ensembles de 2-cohomologie non ab\'elienne correspondants). On note $\eta_X^{\ab}$ l'image de $\eta_X$ dans $\H^2(K,H^{\ab})$.

Sur $\ol{K}$, on a $\ol{X} = \ol{G} / \ol{H}$. Au vu du diagramme
	$$\xymatrix{
		\ol{G} \ar[rd]^{[\ol{H},\ol{H}]} \ar[dd]_{\ol{H}} \\
		& \ol{Z} \ar[ld]^{\ol{H}^{\ab}} \\
		\ol{X},
	}$$
on peut consid\'erer la classe $[\ol{G}] \in \H^1_\et(\ol{X},\ol{H})$ (resp. $[\ol{Z}] \in \H^1_\et(\ol{X},\ol{H}^{\ab})$) du $\ol{H}$-torseur $\ol{G}$ (resp. du $\ol{H}^{\ab}$-torseur $\ol{Z}$) sur $\ol{X}$. Leur {\it type} est par d\'efinition le morphisme
	$$\lambda: H' \to \Pic \ol{X} = \H^1_{\et}(\ol{X},\Gbb_m), \qquad \chi \mapsto \chi_\ast[\ol{G}] = \chi_\ast[\ol{Z}].$$
Comme $\Pic \ol{G} = 0$ ($G$ \'etant simplement connexe) et $\ol{K}[X]^\times = \ol{K}^\times$, il est bien connu que $\lambda$ est un isomorphisme de $\Gamma_K$-modules (voir par exemple \cite[Lemma 6.2]{Bokun} et \cite[d\'emonstration du Theorem 9.5.1]{Sko}).

\begin{lemm} \label{lemDifferential}
	La classe $[\ol{Z}] \in \H^1_\et(\ol{X},\ol{H}^{\ab})$ est $\Gamma_K$-invariante. De plus, notant 
		$$d_2^{0,1}: \H^0(K,\H^1_\et(\ol{X},\ol{H}^{\ab})) \to \H^2(K,H^{\ab})$$
	la diff\'erentielle  tir\'ee de la suite spctrale de Hochschild-Serre 
	\begin{equation} \label{eqDifferentialHSSpec}
		E_2^{pq} = \H^p(K,\H^q_\et(\ol{X},\ol{H}^{\ab})) \Rightarrow \H^{p+q}_\et(X,H^{\ab}),
	\end{equation} 
	on a $d_2^{0,1}([\ol{Z}]) = \eta_X^{\ab}$.
\end{lemm}
\begin{proof}
	Soit $\pi: X \to \Spec K$ le morphisme structural. On consid\`ere la suite spectrale de Leray 
	\begin{equation} \label{eqDifferentialLeraySpec}
		\tilde{E}_2^{pq} = \Ext_{K-\grp}^p(H',\R^q \pi_\ast \Gbb_m) \Rightarrow \H^{p+q}_\et(X,H^{\ab}).
	\end{equation}
	Comme $\ol{K}[X]^\times = \ol{K}^\times$ est divisible, il y a un r\'esultat de Skorobogatov \cite[Proposition 2.3.11]{Sko} qui compare \eqref{eqDifferentialHSSpec} et \eqref{eqDifferentialLeraySpec}. En particulier, on dispose d'un diagramme commutatif dont les fl\`eches verticales sont des isomorphismes: 
	\begin{equation} \label{eqDifferentialDiagram}
		\xymatrix{
			\H^0(K,\H^1_\et(\ol{X},\ol{H}^{\ab})) \ar[d]^{\type} \ar[rr]^{d_2^{0,1}} && \H^2(K,H^{\ab}) \ar@{=}[d] \\
			\Hom(H',\Pic \ol{X})^{\Gamma_K}
			\ar[rr]^{\tilde{d}_2^{0,1}} && \H^2(K,H^{\ab}).
		}
	\end{equation}
	Dans \eqref{eqDifferentialDiagram}, $\tilde{d}_2^{0,1}$ is la  diff\'erentielle tir\'ee de \eqref{eqDifferentialLeraySpec}, et l'isomorphisme $\type: \H^1_\et(\ol{X},\ol{H}^{\ab}) \to \Hom(H',\Pic\ol{X})$ associe \`a chaque classe d'isomorphie de $\ol{H}^{\ab}$-torseurs sur $\ol{X}$ leur type. Comme $\lambda \in \Hom(H',\Pic \ol{X})^{\Gamma_K}$ provient de $[\ol{Z}] \in \H^1_\et(\ol{X},\ol{H}^{\ab})$, on a $[\ol{Z}] \in \H^0(K,\H^1_\et(\ol{X},\ol{H}^{\ab}))$. Il reste donc \`a montrer que $\tilde{d}_2^{0,1}(\lambda) = \eta_X^{\ab}$. 
	
	Maintenant, on suit la preuve de \cite[Theorem 9.5.1]{Sko}. La classe $e(X):=\tilde{d}_2^{pq}(\lambda) \in \H^2(k,H^{\ab})$ est appel\'ee {\it obstruction el\'ementaire}  de $X$, c'est une obstruction \`a l'existence de $H^{\ab}$-torseurs sur $X$ de type $\lambda$. Elle est represent\'ee par le gerbe $\Gcal_\lambda$ des torseurs over $X$ of type $\lambda$, {\it i.e.} pour toute extension finie s\'eparable $L/K$, la fibre $\Gcal_\lambda(L)$ est le groupo\"ide des $H^{\ab}_L$-torseurs sur $X_L$ de type $\lambda$. D'ailleurs, la classe de Springer $\eta_X \in \H^2(K,L_X)$ est represent\'ee par le gerbe $\Gcal_X$, dont la fibre $\Gcal_X(L)$ est pour toute extension finie s\'eparable $L/K$ le groupo\"ide des espaces homog\`enes principaux de $G_L$ dominant $X_L$. Maintenant, soit $Y \to X_L$ est un morphisme $G_L$-equivariant, o\`u $L/K$ est une extension finie s\'eparable et $Y$ est un espace homog\`ene principal de $G_L$. Alors le groupe alg\'ebrique $\Aut_{G_L}(Y/X_L)$ est une $L$-forme de $\ol{H}$, qu'on va noter $H_L$, et $Y$ est un $H_L$-torseur sur $X_L$. Le produit contract\'e $W:=Y \times^{H_L} H_L^{\ab}$ est un $H_L^{\ab}$-torseur sur $X_L$. Comme $\ol{Y} = \ol{G}$, on a $\ol{W} = \ol{G}/[\ol{H},\ol{H}] = \ol{Z}$, qui est un $\ol{H}^{\ab}$-torseur sur $\ol{X}$ de type $\lambda$. La construction $Y \mapsto W$ d\'efinit un morphisme $\Gcal_X \to \Gcal_\lambda$ de $K$-gerbes alg\'ebriques, donc l'application $\H^2(K,L_X) \to \H^2(K,H^{\ab})$ envoie $\eta_X$ sur $e(X)$, {\it i.e.} $d_2^{0,1}([\ol{Z}]) = \tilde{d}_2^{0,1}(\lambda) = e(X) = \eta_X^{\ab}$.
\end{proof}

\begin{prop} \label{propDifferential}
	Soient $G$, $X$ et $\ol{H}$ comme ci-dessus. On utilise le type $\lambda$ de $X$ pour identifier $H'$ \`a  $\Pic \ol{X}$. Alors pour tout $p \ge 0$, la diff\'erentielle $d^{p,1}_2: \H^p(K,H') \to \H^{p+2}(K,\Gbb_m)$	de la suite spectrale de Hochschild-Serre $E_2^{pq} = \H^p(K,\H^q_\et(\ol{X},\Gbb_m)) \Rightarrow \H^{p+q}_\et(X,\Gbb_m)$ est donn\'ee par les cup-produits avec $\eta_X^{\ab} \in \H^2(K,H^{\ab})$.
\end{prop}
\begin{proof}
	Ce calcul est \'essentiellement pris de \cite[Theorem 2.4.4]{Neukirch}. Pour tout faisceau $\Acal$ de groupes ab\'eliens sur $X_{\et}$, on utilisera abusivement la m\^eme notation pour sa restriction \`a $\ol{X}_{\et}$. Soit  $0 \to \Acal \to \Gcal^0_\Acal \to \Gcal^1_\Acal \to \cdots$ la r\'esolution de Godement de $\Acal$, alors la suite spectrale de Hochschild-Serre
	\begin{equation} \label{eqDifferentialHochschildSerre}
		E_2^{pq} = \H^p(K,\H^q_\et(\ol{X},\Acal)) \Rightarrow \H^{p+q}_\et(X,\Acal)
	\end{equation}
	est induite par le complexe double $(\C^\bullet(K,\Gcal^\bullet_\Acal(\ol{X})))$. 
	Sa diff\'erentielle $d_2^{p,1}$ de est obtenue de la partie
	$$\xymatrix{
		\C^p(K,\Gcal_\Acal^1(\ol{X})) \ar[r] & \C^{p+1}(K,\Gcal_\Acal^1(\ol{X})) \\
		& \C^{p+1}(K,\Gcal_\Acal^0(\ol{X})) \ar[u] \ar[r] & \C^{p+2}(K,\Gcal_\Acal^0(\ol{X}))
	}$$
	de ce complexe comme suit. Consid\'erons les suites exactes courtes
	\begin{equation} \label{eqDifferentialExact1}
		0 \to M_\Acal \to N_\Acal \to \H^1_\et(\ol{X},\Acal) \to 0 
	\end{equation}
	et
	\begin{equation} \label{eqDifferentialExact2}
		0 \to \Acal(\ol{X}) \to \Gcal_\Acal^0(\ol{X}) \to M_\Acal \to 0,
	\end{equation}
	de $\Gamma_K$-modules, o\`u $M_\Acal := \Img(\Gcal_\Acal^0(\ol{X}) \to \Gcal_\Acal^1(\ol{X}))$ et $N_\Acal := \Ker(\Gcal_\Acal^1(\ol{X}) \to \Gcal_\Acal^2(\ol{X}))$. Alors $d_2^{p,1}: \H^p(K,\H^1_\et(\ol{X},\Acal)) \to \H^{p+2}(K,\Acal(\ol{X}))$ est la compos\'ee des morphismes connectants
		$$\H^p(K,\H^1_\et(\ol{X},\Acal)) \xrightarrow{\del} \H^{p+1}(K,M_\Acal) \xrightarrow{\delta} \H^{p+2}(K,\Acal(\ol{X}))$$
	induits par \eqref{eqDifferentialExact1} et par \eqref{eqDifferentialExact2}. On prend pour $\Acal$ respectivement $H^{\ab}_X$ et $\Gbb_m$. Alors pour tout caract\`ere $\chi \in H'$, on a un diagramme commutatif
		$$\xymatrix{
			0 \ar[r] & \ol{H}^{\ab} \ar[d]^{\chi} \ar[r] & \Gcal^0_{H^{\ab}_X}(\ol{X})  \ar[d]^{\Gcal^0_\chi} \ar[r] & \Gcal^1_{H^{\ab}_X}(\ol{X})  \ar[d]^{\Gcal^1_\chi} \ar[r] & \cdots \\		
			0 \ar[r] & \ol{K}^\times \ar[r] & \Gcal^0_{\Gbb_m}(\ol{X}) \ar[r] & \Gcal^1_{\Gbb_m}(\ol{X}) \ar[r]	& \cdots 
		}$$
	Le fait que le foncteur de Godement est $\Gamma_K$-equivariant nous permet de d\'efinir des accouplements
		$$\Gcal^p_{H_X^{\ab}}(\ol{X}) \times H' \to \Gcal^p_{\Gbb_m}(\ol{X}), \qquad (u,\chi) \mapsto \Gcal^p_{\chi}(u)$$
	pour tout $p \ge 0$. D'o\`u des accouplements
		$$M_{H_X^{\ab}} \times H' \to M_{\Gbb_m}, \qquad N_{H_X^{\ab}} \times H' \to N_{\Gbb_m}, \qquad \text{et} \qquad \Gcal^0_{H^{\ab}_X}(\ol{X}) \times H' \to \Gcal^0_{\Gbb_m}(\ol{X}),$$
	qui, au vu de \eqref{eqDifferentialExact1} and \eqref{eqDifferentialExact2}, sont compatibles avec les accouplements
		$$\H^1_\et(\ol{X},\ol{H}^{\ab}) \times H' \to \Pic\ol{X}, \qquad (u,\chi) \mapsto \chi_\ast u$$
	et 
		$$\ol{H}^{\ab} \times H' \to \ol{K}^\times, \qquad (h,\chi) \mapsto \chi(h).$$
	Par compatiblit\'e des cup-produits avec les morphismes connectants \cite[Proposition 1.4.3]{Neukirch}, on dispose d'un diagramme commutatif
	\begin{equation} \label{eqDifferentialDiagram2}
		\xymatrix{
			\H^0(K,\H^1_\et(\ol{X},\ol{H}^{\ab})) \ar[d]^{\del} &\times& \H^p(K,H') \ar@{=}[d] \ar[r]^-{\cup} & \H^p(K,\Pic\ol{X})	\ar[d]^{\del} \\
			\H^1(K,M_{H^{\ab}_X}) \ar[d]^{\delta} &\times& \H^p(K,H') \ar@{=}[d] \ar[r]^-{\cup} & \H^{p+1}(K,M_{\Gbb_m}) \ar[d]^{\delta}	\\
			\H^2(K,H^{\ab}) &\times& \H^p(K,H') \ar[r]^-{\cup} & \H^{p+2}(K,\Gbb_m)			
		}
	\end{equation}
	Regardons la classe $[\ol{Z}] \in \H^0(K,\H^1_\et(\ol{X},\ol{H}^{\ab}))$ et soit $y \in \H^p(K,H')$. Par d\'efinition de $\lambda$, on a $[\ol{Z}] \cup y = \lambda_\ast y \in \H^p(K,\Pic \ol{X})$. De l'autre c\^ot\'e, $\delta(\del([\ol{Z}])) = d_2^{0,1}([\ol{Z}]) = \eta_X^{\ab}$ par le lemme \ref{lemDifferential}, d'o\`u \eqref{eqDifferentialDiagram2} donne
	$$d_2^{p,1}(\lambda_\ast y) = \delta(\del(\lambda_\ast y)) = \delta(\del([\ol{Z}] \cup y)) = \delta(\del([\ol{Z}])) \cup y =  \eta_X^{\ab} \cup y,$$
	qui est ce qu'on veut.
\end{proof}

\begin{coro} \label{corDifferential}
	Soit $n \ge 2$ un entier et soit $K$ un corps de caract\'eristique ne divisant pas $n$. Si $K$ contient $\mu_n$ et si $\H^3(K,\Gbb_m)[n] \neq 0$, alors il existe un $K$-espace homog\`ene $X$ de $\SL_m$ tel que la suite exacte \ref{eqExactSequenceX} identifie $\Br_1 X / \Br K$ \`a un sous-groupe stricte de $\H^1(K,\Pic \ol{X})$.
\end{coro}
\begin{proof}
	Par le th\'eor\`eme de Rost-Voedvodsky \cite[Theorem 6.16]{Voedvodsky}, on a un isomorphisme $K^M_\bullet(K)/n \simeq \H^\bullet(K,\mu_n^{\bullet})$ d'anneaux anti-commutatifs gradu\'es, o\`u $K^M_\bullet$ d\'esigne la $K$-th\'eorie de Milnor. D'o\`u $K^M_\bullet(K)/n \simeq \H^\bullet(K,\Zbb/n)$ puisque $K$ contient $\mu_n$. Les g\'en\'erateurs de $K^M_\bullet(K)$ se trouvant en degr\'e $1$, tout \'el\'ement de $\H^r(K,\Zbb/n)$ (o\`u $r \ge 1$) est une somme de symboles $a_1 \cup \cdots \cup a_r$, o\`u $a_1,\ldots,a_r \in \H^1(K,\Zbb/n)$. Soit $c \in \H^3(K,\Gbb_m)[n]$ non nul. Comme $\H^3(K,\mu_n)\to \H^3(K,\Gbb_m)[n]$ est surjectif, $c$ se rel\`eve en un \'el\'ement de $\H^3(K,\mu_n)$, il est donc une somme de symboles de la forme $a \cup b$, o\`u $a \in \H^1(K,\Zbb/n)$ et $b \in \H^2(K,\mu_n)$. Comme $c \neq 0$, l'un de ces symboles est non nul, {\it i.e.} il existe $a \in \H^1(K,\Zbb/n)$ et $b \in \H^2(K,\mu_n)$ tel que $a \cup b \neq 0 \in \H^3(K,\Gbb_m)$. La construction de Demarche-Lucchini Arteche \cite[Corollaire 3.5]{demarche2019reduction} donne un espace homog\`ene $X$ de $\SL_m$ de lien de Springer $L_X = \lien(\mu_n)$ et de classe de Springer $\eta_X = b \in \H^2(K,\mu_n)$. Au vu de la proposition \ref{propDifferential}, on a $\Pic \ol{X} = \Zbb/n$, et l'image de $a \in \H^1(K,\Zbb/n)$ dans $\H^3(K,\Gbb_m)$ vaut $a \cup b \neq 0$. Par exactitude de \eqref{eqExactSequenceX}, $a$ ne provient pas de $\Br_1 X$.
\end{proof}

\begin{exam}
	Pour $K = \Cbb\ps{t}\ps{x}\ps{y}$, on a $\H^3(K,\Gbb_m) = \Qbb/\Zbb$ (voir le \S \ref{section2} ci-dessous), d'o\`u un exemple num\'erique du corollaire \ref{corDifferential}.
\end{exam}

\section{Le contre-exemple} \label{section2}

Notre principal r\'esultat est le

\begin{thm} \label{thmMain}
	Il existe un espace homog\`ene $X$ de $\SL_m$ sur $K := \Cbb\ps{t} \ps{x} \ps{y}$, \`a stabilisateur g\'eom\'etrique fini tel que la suite exacte \eqref{eqExactSequenceXc} identifie $\Br_1 X^c / \Br K$ \`a un sous-groupe stricte de $\Sha^1_\omega(K,\Pic\ol{X})$.
\end{thm} 

Mentionnons par passage quelques th\'eor\`emes de dualit\'es pour des corps locaux sup\'erieurs. Les corps 0-locaux sont par d\'efinition les corps finis et le corps $\Cbb\ps{t}$. Pour $d \ge 1$, on appelle corps $d$-local tout corps complet pour une valuation discr\`ete de corps r\'esiduel un corps $(d-1)$-local. Si $K$ est un corps $d$-local, on notera $K_d := K$, et $K_{i-1}$ le corps r\'esiduel de $K_i$ pour $i \in \{1,\ldots,d\}$. Les points 1. et 2. de la proposition \ref{propLocalDuality} ci-dessous se demontrent de mani\`ere similaire \`a celle comme dans \cite{milne2006duality}[Theorem 2.17, Lemma 2.18], compte tenu du fait que $\mu \simeq \Qbb/\Zbb$ sur toutes les extensions de $\Cbb$, et que le  module galoisien $\Gbb_m/\mu$ est uniquement divisible (donc {\it cohomologiquement trivial}). On pourra consulter \cite[Th\'eor\`eme 8.9]{harari2017cohomologie} pour le point 3. Le point 4. se trouve dans \cite[Proposition 3.5]{diego}.

\begin{prop} \label{propLocalDuality}
	Soit $K$ un corps $d$-local $(d \ge 1$) avec $K_0 = \Cbb\ps{t}$.
	\begin{enumerate}
		\item On dispose d'un isomorphisme (non canonique) $\H^{d+1}(K,\Gbb_m) \simeq \Qbb/\Zbb$. 
		
		\item Pour tout $\Gamma_K$-module fini $M$ et $r \in \{0,1,\ldots,d+1\}$, le cup-produit
		    $$\H^r(K,M') \times \H^{d+1-r}(K,M) \to \Qbb/\Zbb$$
		est une dualit\'e parfaite de groupes finis.
		
		\item Pour toute extension finie $L/K$, la restriction $\H^{d+1}(K,\Gbb_m) \to \H^{d+1}(L,\Gbb_m)$ est la multiplication par $[L:K]$ de $\Qbb/\Zbb$.
		
		\item Supposons $T = 2$. Soit $T$ un $K$-tore et $S$ son tore dual, c'est-\`a-dire $\Hom(T,\Gbb_m) = \Hom(\Gbb_m,S)$. Il existe un accouplement parfait
		    $$\varprojlim_n\H^0(K,S) / n \times \H^2(K,T) \to \Qbb/\Zbb.$$
		o\`u $\varprojlim_n\H^0(K,S) / n$ est profini et $\H^2(K,T)$ est discrete de torsion.
	\end{enumerate}
\end{prop}

\begin{coro} \label{corH3KMuN}
	Soit $K = \Cbb\ps{t}\ps{x}\ps{y}$. Pour tout entier $n$, $\H^3(K,\mu_n)$ s'identifie au sous-groupe $\frac{1}{n}\Zbb/\Zbb$ de $\H^3(K,\Gbb_m) \simeq \Qbb/\Zbb$.
\end{coro}
\begin{proof}
	Le point 4. de la proposition \ref{propLocalDuality} appliqu\'e au tore $T = \Gbb_m$ donne une dualit\'e parfaite entre $\Br K$ et $\varprojlim_n K^\times/K^{\times n}$. Rappelons que pour tout corps $k$, on a $k\ps{t}^\times = \Zbb \times k^\times \times (1 + tk\pssq{t})$, et $(1 + tk\pssq{t})$ est divisible lorsque $k$ est de caract\'eristique nulle. Donc $K^\times = \Zbb^3 \times \Cbb^\times \times D$ avec $D$ divisible. Comme $\Cbb^\times$ est aussi divisible, $\varprojlim_n K^\times/K^{\times n} = \hat{\Zbb}^3$, d'o\`u $\Br K = (\Qbb/\Zbb)^3$.
	
    Pour tout entier $n$, la suite de Kummerdonne lieu \`a une suite exacte
		$$0 \to \Br K / n \to \H^3(K,\mu_n) \to \H^3(K,\Gbb_m)[n] \to 0.$$
	Or $\Br K = (\Qbb/\Zbb)^3$ est divisible, d'o\`u $\H^3(K,\mu_n) \simeq \H^3(K,\Gbb_m)[n]$. Finalement, $\H^3(K,\Gbb_m) \simeq \Qbb/\Zbb$ au vu de la proposition \ref{propLocalDuality}, d'o\`u le r\'esultat.
\end{proof}

Pour d\'emontrer le th\'eor\`eme \ref{thmMain}, la premi\`ere \'etape est de construire un $\Gamma_K$-module $H$ avec $\Sha^1_\omega(K,H') \neq 0$. La construction suivante est inspir\'ee par un exemple de Serre \cite[Chapitre III, 4.7, Lemme 7]{serre2013cohomologie}, et est \'egalemente r\'eapparue par exemple dans des articles de Borovoi-Kunyavskii \cite{Bokun}, Demarche-Lucchini Arteche-Neftin \cite[Lemma 5.5]{DLN} et Rivera-Mesas \cite[Lemma 4.1]{RM}.

\begin{prop} \label{propData}
	Soient $n$ un entier, $K$ un corps contenant $\mu_n$ et de caract\'eristique ne divisant pas $n$, et $L/K$ une extension finie galoisienne. Notons $\gfrak = \Gal(L/K)$, $N = \PGCD(n,|\gfrak|)$ et $N' = \PGCD(n,\exp(\gfrak))$. Soit $j: \mu_n \hookrightarrow \R_{L/K} \mu_n$ l'inclusion canonique et soit $H$ le $\Gamma_K$-module qui remplit la suite exacte
	\begin{equation} \label{eqExactSequenceH}
		1 \to \mu_n \xrightarrow{j} \R_{L/K} \mu_n \to H \to 1.
	\end{equation}
	Alors $\Sha^1_\omega(K,H')$ est un sous-groupe cyclique d'ordre $N/N'$ dont un g\'en\'erateur est $\delta'(\exp(\gfrak))$, o\`u $\delta': \Zbb/n \to \H^1(K,H')$ est le morphisme connectant induit par la suite exacte duale \`a \eqref{eqExactSequenceH}:
	\begin{equation} \label{eqExactSequenceHDual}
		1 \to H' \to (\Zbb/n)[\gfrak] \xrightarrow{j'} \Zbb/n \to 1.
	\end{equation}
\end{prop}
\begin{proof}
	On v\'erifie sans peine que $j'$ est l'application d'augmentation, {\it i.e.} si $(e_g)_{g \in \gfrak}$ d\'esigne la $(\Zbb/n)$-base canonique de $(\Zbb/n)[\gfrak]$, alors $j'(e_g) = 1$ pour tout $g \in \gfrak$. 
	
	Au vu de \eqref{eqSha1Omega}, $\Sha^1_\omega(K,H') = \Sha^1_\omega(\gfrak, H') = \Ker\tuple{\H^1(\gfrak,A) \to \prod_{g \in \gfrak} \H^1(\pair{g},A)}$. Pour tout $g \in \gfrak$, \eqref{eqExactSequenceHDual} donne un diagramme commutatif \`a lignes exactes:
	\begin{equation} \label{eqExactSequenceSha1Omega}
		\xymatrix{
			\H^0(\gfrak, (\Zbb/n)[\gfrak]) \ar[r]^-{j'} \ar[d] & \Zbb/n \ar@{=}[d] \ar[r]^-{\delta'} & \H^1(\gfrak,H') \ar[d] \ar[r] & 0 \\
			\H^0(\pair{g}, (\Zbb/n)[\gfrak]) \ar[r]^-{j'} & \Zbb/n \ar[r]^-{\delta'} & \H^1(\pair{g},H') \ar[r] & 0,
		}
	\end{equation}
	o\`u $\H^1(\gfrak, (\Zbb/n)[\gfrak]) = \H^1(\pair{g}, (\Zbb/n)[\gfrak]) = 0$ par le lemme de Shapiro. En particulier, tout \'el\'ement  $c \in \H^1(\gfrak,H')$ s'\'ecrit sous la forme $\delta'(m)$, o\`u $m \in \Zbb/n$. Comme $j'$ est l'application d'augmentation et comme la $(\Zbb/n)$-base $(e_g)_{g \in \gfrak}$ de $(\Zbb)[\gfrak]$ est permut\'ee par $\gfrak$, on voit que l'image de $\H^0(\gfrak, (\Zbb/n)[\gfrak])$ par $j'$ est $|\gfrak|(\Zbb/n)$, et celle de $\H^0(\pair{g}, (\Zbb/n)[\gfrak])$ est $\ord(g)(\Zbb/n)$. Donc, $c \in \Sha^1_\omega(\gfrak,H')$ si et seulement si $m \in \ord(g)(\Zbb/n)$ pour tout $g \in \gfrak$, {\it i.e.} $m \in \exp(\gfrak)(\Zbb/n) = N'\Zbb/n$. D'autre part, $c = 0$ si et seulement si $m \in |\gfrak|(\Zbb/n) = N(\Zbb/n)$. D'o\`u $\Sha^1_\omega(\gfrak,H') \simeq \frac{N'(\Zbb/n)}{N(\Zbb/n)}$ est un groupe cyclique d'ordre $N/N'$, engendr\'e par $\delta'(\exp(\gfrak))$.
\end{proof}

\begin{lemm} \label{lemShapiroRestriction}
	Avec les donn\'ees de la proposition \ref{propData}, on a une suite exacte longue	
		$$\cdots \to \H^{r-1}(K,H) \xrightarrow{\delta} \H^r(K,\mu_n) \xrightarrow{\res} \H^r(L,\mu_n) \to \H^r(K,H) \xrightarrow{\delta} \H^{r+1}(K,\mu_n) \to \cdots,$$
	o\`u $\delta$ sont les morphismes connectants induits par \eqref{eqExactSequenceH}.
\end{lemm}
\begin{proof}
    C'est la suite exacte longue induite par \eqref{eqExactSequenceH}, compte tenue du fait que $\H^r(K,\R_{L/K} \mu_n) = \H^r(L,\mu_n)$ pour tout $r \ge 1$ par le lemme de Shapiro, et de la compatibitli\'e de $j'_\ast: \H^r(K,\mu_n) \to \H^r(K,\R_{L/K}\mu_n)$ avec la restriction (voir \cite[Proposition 1.6.5]{Neukirch}).
\end{proof}

\begin{proof} [D\'emosntration du th\'eor\`eme \ref{thmMain}]
	On consid\`ere  les donn\'ees suivantes: $K = \Cbb\ps{t}\ps{t_1}\ps{t_2}$, $L = K(\sqrt{t_1},\sqrt{t_2})$ et $n = 4$. Avec les notations de la proposition \ref{propData}, on a $N' = 2$ et $N = 4$. Regardons l'\'element $a = \frac{1}{4} \in \Qbb/\Zbb \simeq \H^3(K,\Gbb_m)$. Sa restriction \`a $\H^3(L,\Gbb_m)$ vaut $0$ au vu de la proposition \ref{propLocalDuality}. En vertu du corollaire \ref{corH3KMuN}, on a $a \in \H^3(K,\mu_4)$, et $\res(a)  = 0 \in \H^3(L,\mu_4)$. Par le lemme \ref{lemShapiroRestriction}, il existe $x \in \H^2(K,H)$ tel que $a = \delta(x)$.
	
	On prend $y = \delta'(2) \in \H^1(K,H')$. Alors $y \in \Sha^1_\omega(K,H')$, et $x \cup y = \pm 2a \in \H^3(K,\Gbb_m)$ par compatibilt\'e des cup produits avec les morphismes connectants. Or $2a \neq 0$, donc $x \cup y \neq 0$. Soit maintenant $X$ un espace homog\`ene de $\SL_m$ \`a de lien de Springer $\lien(H)$ et de classe de Springer $\eta_X = x$, d'apr\`es la construction de Demarche-Lucchini Arteche \cite[Corollaire 3.5]{demarche2019reduction}. Par la proposition \ref{propDifferential}, $\Pic \ol{X} = H'$, et l'image de $y \in \H^1(K,H')$ dans $\H^3(K,\Gbb_m)$ est $x \cup y \neq 0$, donc l'exactitude de \eqref{eqExactSequenceXc} assure que $y$ ne provient pas de $\Br_1 X^c$.
\end{proof}

	\bibliographystyle{alpha-fr}
	\selectlanguage{french}
	\bibliography{ref}
	
\end{document}